# An Algorithm for Verifying Some Norm Identities in Inner-Product Spaces

Muneerah Al Nuwairan

Department of mathematics and statistics

King Faisal university, P.O. Box 400

Al ahssa 31982, Saudi Arabia.

E-mail: msalnuwairan@kfu.edu.sa

**Abstract:** In this paper, we provide an algorithm for verifying the validity of identities of the form $\sum_{A \subseteq \bar{n}} c_A \|x_A\|^2 = 0$, where $x_A = \sum_{i \in A} x_i$ and $\bar{n} = \{1, \ldots, n\}$ in inner-product spaces. Such algorithm is used to verify the validity, in inner-product spaces, for a number of identities. These include a generalization of the parallelopiped law. We also show that such identities hold only in inner-product spaces. Thus, the algorithm can be used to deduce characterizations of inner-product spaces.

## 1. Introduction:

Throughout this paper, let $I$ be an index set, and $\{x_i : i \in I\}$ be a subset of elements of a vector space $\mathcal{H}$. For $A \subseteq I$ denote by $x_A$ the sum of vectors $x_i, i \in A$. i.e. $x_A = \sum_{i \in A} x_i$. The notation $\bar{n}$ is used to denote the set $\{1, 2, \ldots, n\}$. For a finite set $A$, we use $|A|$ to denote the cardinality of $A$, and the standard notation for binomial coefficients, $\binom{n}{k} = \frac{n!}{k!(n-k)!}$ is used.

The algorithm given in Theorem 2.3 is meant to test the validity of identities of the form $\sum_{A \subseteq \bar{n}} c_A \|x_A\|^2 = 0$ in inner-product spaces by converting the verification of such an identity to verifying numerical equalities. The algorithm is illustrated in Section 3 by using it to derive several identities. Notable among these results is a generalization of the parallelopiped law, which is deduced, in Corollary 3.5, from a more general result. In Section 4, we prove that all the identities that can be verified by this algorithm only hold in inner-product spaces. Thus, the algorithm can be used to derive characterizations of norms defined by an inner product. This is the application chosen for discussion in the paper.

Investigating norm identities that are satisfied only by norms induced by inner products dates back to the late 19[th] century (see [2] Introduction). Fréchet [5] showed that a normed space is an inner product space if and only if

$$\|x + y + z\|^2 - \|x + y\|^2 - \|x + z\|^2 - \|y + z\|^2 + \|x\|^2 + \|y\|^2 + \|z\|^2 = 0 \quad (1)$$

Jordan and von Neumann [7] showed that the norm is induced by an inner product if and only if the parallelogram law

$$\|x + y\|^2 + \|x - y\|^2 = 2\|x\|^2 + 2\|y\|^2 \quad (2)$$

holds for all $x, y$. For a proof of Jordan and von Neumann's result see [6, Thm 4.3.6]. The study of characterizations of inner-product spaces continue to be an active field (see e.g. [1],[3],[4], and [8]). The author hopes that the algorithm presented here will lead to new characterization like the ones given in Corollary 4.5.

This paper is self-contained, the only results needed are the Jordan and von Neumann characterization mentioned above and a Lemma due to Fréchet on finite difference given below (see e.g. [2, Lemma 1.1, Lemma 1.2] or [9 Lemma 1]). For completeness, we include a proof.

**Lemma 1.1:**

1. Let $\{a_m\}_{m \in \mathbb{Z}}$ be a sequence of real numbers and $n \in \mathbb{N}$. If for each $l \in \mathbb{Z}$,
$$\sum_{k=0}^{n} \binom{n}{k} (-1)^{n-k} a_{l+k} = 0$$
then there is a polynomial $P$ of degree less than $n$ such that $a_l = P(l)$ for all $l \in \mathbb{Z}$.
2. Let $g: \mathbb{R} \to \mathbb{R}$ be continuous map and $n \in \mathbb{N}$. If for all $r, s \in \mathbb{R}$
$$\sum_{k=0}^{n} \binom{n}{k} (-1)^{n-k} g(r + ks) = 0$$
then $g$ is a polynomial of degree less than $n$.

**Proof:**

1. Let $E$ be the operator on sequences defined by $(Ea)_n = a_{n+1}$, so $(E^k a)_n = a_{n+k}$. Let $I$ be the identity operator on sequences, then the hypothesis can be written as $(E - I)^n a = 0$. We show that this implies $a_l = P(l)$ for a polynomial of degree less than $n$ by induction on $n$. For $n = 1$ our hypothesis is $a_{l+1} - a_l = 0, \forall l \in \mathbb{Z}$, so the sequence is constant $a_l = c \quad \forall l \in \mathbb{Z}$. By taking $P \equiv c$ of degree 0 the result follows in this case. Suppose the result is true for $n \geq 1$ and that $(E - I)^{n+1} a = 0$.
Let $b = (E - I)a$ then $(E - I)^n b = 0$. By induction hypothesis, there exists a polynomial $Q$ of degree less than $n$ such that $b_l = Q(l) \quad \forall l \in \mathbb{Z}$. Write
$$Q(x) = \sum_{i=0}^{n-1} c_i x^{(i)} \quad \text{where} \quad x^{(i)} = x(x-1)\cdots(x-i+1).$$
By direct computation,
$$(E - I)x^{(i)} = ix^{(i-1)} \text{ for } i \geq 1.$$
Let $P_1(x) = \sum_{i=0}^{n-1} \frac{c_i}{i+1} x^{(i+1)}$. Clearly $(E - I)P_1 = Q$, so
$$a_{l+1} - a_l = b_l = Q(l) = P_1(l + 1) - P_1(l) \quad \text{for all } l \in \mathbb{Z}$$
i.e.
$$a_{l+1} - P_1(l + 1) = a_l - P_1(l).$$
and there is a constant $c_0$ such that $a_l = P_1(l) + c_0, \quad \forall l \in \mathbb{Z}$.
Thus $P(x) = P_1(x) + c_0$ is the desired polynomial of degree less than $n + 1$.
2. For each integer $q \geq 0$, we define a sequence $a^q$ by $a_m^q = g\left(\frac{m}{2^q}\right), m \in \mathbb{Z}$. Note that
$$\sum_{k=0}^{n} \binom{n}{k} (-1)^{n-k} a_{l+k}^q = \sum_{k=0}^{n} \binom{n}{k} (-1)^{n-k} g\left(\frac{l+k}{2^q}\right)$$

$$= \sum_{k=0}^{n} \binom{n}{k} (-1)^{n-k} g(r+ks).$$

where $r = \frac{l}{2^q}$ and $s = \frac{1}{2^q}$.

So, $\sum_{k=0}^{n} \binom{n}{k}(-1)^{n-k} a_{l+k}^q = 0$ for any $l \in \mathbb{Z}$. Thus, for each integer $q \geq 1$, there exist a polynomial $P_q$ of degree less than $n$ such that $P_q(m) = a_m^q = g\left(\frac{m}{2^q}\right)$, for all $m \in \mathbb{Z}$. Since, for all $m \in \mathbb{Z}$,

$$P_0(m) = g(m) = P_q(2^q m)$$

the two polynomials $P_q(2^q x)$ and $P_0(x)$ have the same value at infinitely many points therefore, they are identical. So $P_q(y) = P_0\left(\frac{y}{2^q}\right)$ for all $y$, thus

$$g\left(\frac{m}{2^q}\right) = P_q(m) = P_0\left(\frac{m}{2^q}\right).$$

The two continuous functions $g$ and $P$ coincide on a dense set (the dyadic rationals) hence, must be identical and $g$ is a polynomial of degree less than $n$. ∎

## 2. A Test for A Class of Norm Identities in Inner Product Spaces.

In this section, a test that can be used to verify the validity of certain identities in inner product spaces is provided. Recall the familiar identity for inner product norms

$$\|x+y\|^2 = \|x\|^2 + \|y\|^2 + 2Re\langle x, y \rangle \qquad (3)$$

Or equivalently

$$2\,Re\langle x, y \rangle = \|x+y\|^2 - \|x\|^2 - \|y\|^2 \qquad (4)$$

**Theorem 2.1:** Let $\mathcal{H}$ be an inner product space, and $x_1, x_2, \cdots, x_n$ be elements in $\mathcal{H}$. The equality

$$\|x_1 + x_2 + \cdots + x_n\|^2 = \sum_{|A|=2, A \subseteq \overline{n}} \|x_A\|^2 - (n-2) \sum_{i=1}^{n} \|x_i\|^2 \qquad (5)$$

holds for all positive integer $n \geq 2$.

**Proof:**

By induction on $n$. For $n=2$ the equality is simply $\|x_1 + x_2\|^2 = \|x_1 + x_2\|^2$. Suppose the equality holds for $n$ then, from (3)

$$\|x_1 + x_2 + \cdots + x_{n+1}\|^2 = \|x_1 + x_2 + \cdots + x_n\|^2 + \|x_{n+1}\|^2 + 2Re\langle x_1 + \cdots + x_n, x_{n+1}\rangle \quad (6)$$

But

$$2Re\langle x_1 + \cdots + x_n, x_{n+1}\rangle = \sum_{i=1}^{n} 2Re\langle x_i, x_{n+1}\rangle$$

$$= \sum_{i=1}^{n} (\|x_i + x_{n+1}\|^2 - \|x_i\|^2 - \|x_{n+1}\|^2) \qquad (7)$$

$$= \sum_{|A|=2, A \subseteq \overline{n+1}, n+1 \in A} \|x_A\|^2 - \sum_{i=1}^{n} \|x_i\|^2 - n\|x_{n+1}\|^2$$

By induction hypothesis, we have

$$\|x_1 + x_2 + \cdots + x_n\|^2 = \sum_{|A|=2, A\subseteq \overline{n}}\|x_A\|^2 - (n-2)\sum_{i=1}^{n}\|x_i\|^2 \qquad (8)$$

Substituting (7) and (8) in (6), we obtain

$$\|x_1 + x_2 + \cdots + x_{n+1}\|^2 = \left(\sum_{|A|=2, A\subseteq \overline{n}}\|x_A\|^2 - (n-2)\sum_{i=1}^{n}\|x_i\|^2\right) + \|x_{n+1}\|^2$$

$$+ \left(\sum_{|A|=2, A\subseteq \overline{n+1}, n+1\in A}\|x_A\|^2 - \sum_{i=1}^{n}\|x_i\|^2 - n\|x_{n+1}\|^2\right)$$

$$= \sum_{|A|=2, A\subseteq \overline{n+1}}\|x_A\|^2 - (n-1)\sum_{i=1}^{n+1}\|x_i\|^2 .$$

Thus, the relation (5) is true for $n + 1$. By induction, (5) is true for all $n \in \mathbb{N}, n \geq 2$ ∎

By using Theorem 2.1 to substitute for $\|x_A\|^2$ where $|A| > 2$, the equal expression

$$\sum_{|B|=2, B\subseteq A}\|x_B\|^2 - (|A|-2)\sum_{i\in A}\|x_i\|^2$$

any identity of the form $\sum_{B\subseteq \overline{n}} c_A \|x_A\|^2 = 0$ can be converted to an identity of the form

$$\sum_{B\subseteq \overline{n}, 1\leq |B|\leq 2} a_B \|x_B\|^2 = 0$$

Thus, testing the validity of the identity $\sum_{A\subseteq \overline{n}} c_A \|x_A\|^2 = 0$ is transformed into testing the validity of the equivalent identity $\sum_{B\subseteq \overline{n}, 1\leq |B|\leq 2} a_B \|x_B\|^2 = 0$. The verification of the latter identity can be reduced to verifying that all its coefficients are zero as shown by the following Lemma.

**Lemma 2.2:** Let $\mathcal{H}$ be an inner product space of dimension at least two. The identity $\sum_{B\subseteq \overline{n}, 1\leq |B|\leq 2} a_B \|x_B\|^2 = 0$ holds for all $x_1, \dots, x_n \in \mathcal{H}$ if and only if $a_B = 0$ for all $B \subseteq \overline{n}$ with $1 \leq |B| \leq 2$.

**Proof:**

($\Leftarrow$) Clearly if $a_B = 0$ for all $B \subseteq \overline{n}$ with $1 \leq |B| \leq 2$ then

$$\sum_{B\subseteq \overline{n}, 1\leq |B|\leq 2} a_B \|x_B\|^2 = 0 \text{ for all } x_1, \dots, x_n \in \mathcal{H}$$

($\Rightarrow$) Suppose that $\sum_{B\subseteq \overline{n}, 1\leq |B|\leq 2} a_B \|x_B\|^2 = 0$ for all $x_1, \dots, x_n \in \mathcal{H}$. Pick $u, v \in \mathcal{H}$ orthogonal unit vectors. For $0 \leq \theta \leq 2\pi$ let $u_\theta = \cos \theta\ u + \sin \theta\ v$. For $1\leq i \neq j \leq n$, let $x_i = u$, $x_j = u_\theta$ and $x_k = 0$ for $k \in \overline{n}\setminus\{i,j\}$.

Since $\|u + u_\theta\|^2 = 2 + 2\cos \theta$ we have

$$0 = \sum_{B\subseteq \overline{n}, 1\leq |B|\leq 2} a_B \|x_B\|^2$$

$$= \sum_{k \notin \{i,j\}} a_{\{i,k\}} + \sum_{k \notin \{i,j\}} a_{\{j,k\}} + a_{\{i,j\}}(2 + 2\cos\theta) + a_{\{i\}} + a_{\{j\}}$$

Choose $\theta_1, \theta_2 \in [0, 2\pi]$ with $\cos\theta_1 \neq \cos\theta_2$, subtract the above equality at $\theta = \theta_2$ from the same equality at $\theta = \theta_1$ to obtain $a_{\{i,j\}} = 0$. Since $i \neq j$ were arbitrary elements of $\overline{n}$, we obtain that $a_B = 0$ for all $B \subseteq \overline{n}$ with $|B| = 2$. Using this and our assumption of the validity of the identity, we obtain $0 = \sum_{B \subseteq \overline{n}, 1 \leq |B| \leq 2} a_B \|x_B\|^2 = \sum_{j=1}^{n} a_{\{j\}} \|x_j\|^2$. Finally, for each $1 \leq i \leq n$, take $x_i$ to be a unit vector and $x_j = 0$ for $j \neq i$ in the last identity to get

$$0 = \sum_{B \subseteq \overline{n}, 1 \leq |B| \leq 2} a_B \|x_B\|^2 = \sum_{i=1}^{n} a_{\{j\}} \|x_j\|^2 = a_{\{i\}}.$$

Since $1 \leq i \leq n$ was arbitrary, then $a_B = 0$ for all $B \subseteq \overline{n}$ with $|B| = 1$. ∎

The following Theorem uses Theorem 2.1 and a modified version of Lemma 2.2 to test validity of identities of the form $\sum_{A \subseteq \overline{n}} c_A \|x_A\|^2 = 0$. The modification allows us to avoid the need to compute $a_B$ for $B \subseteq \overline{n}$ with $|B| = 1$ which simplifies the application of the validity test.

**Theorem 2.3:** Let $\mathcal{H}$ be an inner product space of dimension at least two. Given an expression $\sum_{A \subseteq \overline{n}} c_A \|x_A\|^2$, let $\sum_{B \subseteq \overline{n}, 1 \leq |B| \leq 2} a_B \|x_B\|^2$ be the result of replacing $\|x_A\|^2$ where $|A| > 2$ by

$$\sum_{|B|=2, B \subseteq A} \|x_B\|^2 - (|A| - 2) \sum_{i \in A} \|x_i\|^2.$$

The identity $\sum_{B \subseteq \overline{n}} c_B \|x_B\|^2 = 0$ holds if and only if

1) $a_B = 0$ for all $B \subseteq \overline{n}$ with $|B| = 2$, and
2) For each unit vector $u \in \mathcal{H}$, and for each $1 \leq i \leq n$, if $x_i = u$ and $x_j = 0$ for $j \neq i$ then for this choice of $x_k$'s we have $\sum_{A \subseteq \overline{n}} c_A \|x_A\|^2 = 0$.

**Proof:**

($\Rightarrow$) If the identity $\sum_{B \subseteq \overline{n}} c_B \|x_B\|^2 = 0$ holds then $\sum_{B \subseteq \overline{n}, 1 \leq |B| \leq 2} a_B \|x_B\|^2 = 0$ is also an identity, so by Lemma 2.2, we have $a_B = 0$ for all $B \subseteq \overline{n}$ with $1 \leq |B| \leq 2$. i.e. 1) holds. Since $\sum_{B \subseteq \overline{n}} c_B \|x_B\|^2 = 0$ holds for any choice of values of the $x_k$'s, 2) also holds.

($\Leftarrow$) Assume 1) and 2) hold, and $u$ be any unit vector in $\mathcal{H}$, for each $1 \leq i \leq n$, let $x_i = u$ and $x_j = 0$ for $j \neq i$ then from 2)

$$0 = \sum_{B \subseteq \overline{n}} c_B \|x_B\|^2 = \sum_{B \subseteq \overline{n}, 1 \leq |B| \leq 2} a_B \|x_B\|^2 = \sum_{B \subseteq \overline{n}, 1 = |B|} a_B \|x_B\|^2 = a_{\{i\}}$$

where we have used 1) in the second equality. Since $1 \leq i \leq n$ was arbitrary, we have $a_B = 0$ for all $B \subseteq \overline{n}$ with $1 = |B|$. This together with 1) gives us that $a_B = 0$ for all $B \subseteq \overline{n}, 1 \leq |B| \leq 2$. Thus $\sum_{B \subseteq \overline{n}} c_B \|x_B\|^2 = \sum_{B \subseteq \overline{n}, 1 \leq |B| \leq 2} a_B \|x_B\|^2 = 0$ is an identity. ∎

3. **Some Deduced Identities.**

In this section, we provide some example applications of Theorem 2.3. For future reference we have listed the identities as Lemmas rather than examples.

**Lemma 3.1**: Let $\mathcal{H}$ be an inner product space. The identity

$$\binom{n-2}{k-2}\|x_1 + \cdots + x_n\|^2 = \sum_{|A|=k, A\subseteq \overline{n}}\|x_A\|^2 - \binom{n-2}{k-1}\sum_{i=1}^{n}\|x_i\|^2 \quad (9)$$

holds for all $n > k \geq 2$ and all $x_1, \ldots, x_n \in \mathcal{H}$.

**Proof:**

We start by converting the identity in (9) to one with zero on one of the sides to get

$$\binom{n-2}{k-2}\|x_1 + \cdots + x_n\|^2 - \sum_{|A|=k, A\subseteq \overline{n}}\|x_A\|^2 + \binom{n-2}{k-1}\sum_{i=1}^{n}\|x_i\|^2 = 0 \quad (10)$$

Using Theorem 2.1 to replace $\|x_A\|^2$ for each $A \subseteq \overline{n}$ with $|A| = k$ by the equivalent expression $\sum_{|B|=2, B\subseteq A}\|x_B\|^2 - (k-2)\sum_{i\in A}\|x_i\|^2$ and replace $\|x_1 + \cdots + x_n\|^2$ by $\sum_{|B|=2, B\subseteq \overline{n}}\|x_B\|^2 - (n-2)\sum_{i=1}^{n}\|x_i\|^2$ converts the left hand side of (10) to

$$\binom{n-2}{k-2}\left(\sum_{|B|=2, B\subseteq \overline{n}}\|x_B\|^2 - (n-2)\sum_{i=1}^{n}\|x_i\|^2\right) -$$

$$\sum_{|A|=k, A\subseteq \overline{n}}\left(\sum_{|B|=2, B\subseteq A}\|x_B\|^2 - (k-2)\sum_{i\in A}\|x_i\|^2\right) + \binom{n-2}{k-1}\sum_{i=1}^{n}\|x_i\|^2 \quad (11)$$

By writing (11) in the form $\sum_{B\subseteq \overline{n}, 1\leq |B|\leq 2} a_B \|x_B\|^2$, we get for a set $B$ of size 2,

$$a_B = \binom{n-2}{k-2} - \sum_{|A|=k, B\subseteq A\subseteq \overline{n}} 1.$$

The number of $k$-subset $A \subseteq \overline{n}$ containing $B$ is equal to the number of $(k-2)$-subsets $C$ of $\overline{n}\backslash B$ so is $\binom{n-2}{k-2}$. Thus for $B \subseteq \overline{n}$ of cardinality 2,

$$a_B = \binom{n-2}{k-2} - \binom{n-2}{k-2} = 0$$

giving us that condition (1) of Theorem 2.3. To verify condition (2), let $u$ to be a unit vector, for each $1 \leq i \leq n$, let $x_i = u$ and $x_j = 0$ for $i \neq j \in \overline{n}$, then substituting this choice of $x_i$'s in the left hand side of (10) (our original identity), the equation in (10) becomes

$$\binom{n-2}{k-2} - \binom{n-1}{k-1} + \binom{n-2}{k-1} = 0$$

The middle term in the above equation was computed by noting that $\|x_A\|^2 = 1$ if $i \in A$ and 0 otherwise, and that the number of k-subsets $A \subseteq \overline{n}$ containing $\{i\}$ is equal to the number of $(k-1)$-subsets of $\overline{n}\backslash\{i\}$ which is $\binom{n-1}{k-1}$. By a well-known recurrence relation for the binomial coefficients (note $k \geq 2$). Thus condition (2) also holds. ∎

**Lemma 3.2:** Let $\mathcal{H}$ be an inner product space of dimension at least 2. For $n \geq 3$ and $x_1, \ldots, x_n$ in $\mathcal{H}$, we have

$$\sum_{I\subseteq\overline{n}}(-1)^{|I|}\|x_I\|^2 = 0 \qquad (12)$$

**Proof:**

Splitting the sum in (12) into a sum over $|I| \geq 2$ and a second sum over $|I| = 1$, then using Theorem 2.1, to substitute for $\|x_I\|^2$ in the first sum, the LHS of (12) becomes

$$\sum_{I\subseteq\overline{n},|I|\geq 2}(-1)^{|I|}\left(\sum_{|B|=2,B\subseteq I}\|x_B\|^2 - (|I|-2)\sum_{i\in I}\|x_i\|^2\right) + \sum_{I\subseteq\overline{n},|I|=1}(-1)\|x_I\|^2 \qquad (13)$$

Let us write (13) in the form

$$\sum_{B\subseteq\overline{n},1\leq |B|\leq 2} a_B \|x_B\|^2$$

For a set $B$ with $|B| = 2$, we have $a_B = \sum_{I\supseteq B}(-1)^{|I|}$. As the number of sets $I \subseteq \overline{n}$ of cardinality $k$ containing $B$ equals the number of ways of choosing the $k-2$ elements of $I \setminus B$ from the $n-2$ elements of $\overline{n}\setminus B$, this sum is

$$a_B = \sum_{k=2}^{n}(-1)^k \binom{n-2}{k-2} = \sum_{l=0}^{n-2}(-1)^l \binom{n-2}{l} = (1-1)^{n-2} = 0$$

So, condition 1) of Theorem 2.3 is satisfied. To verify condition (2), let $u$ be a unit vector, for each $1 \leq i \leq n$ let $x_i = u$ and $x_j = 0$ for $j \neq i$. By substituting this choice in (12), as the number of sets $I \subseteq \overline{n}$ containing $i$ and of cardinality $k$ is $\binom{n-1}{k-1}$, we have

$$\sum_{I\subseteq\overline{n}}(-1)^{|I|}\|x_I\|^2 = \sum_{I\subseteq\overline{n},i\in I}(-1)^{|I|}1 = \sum_{k=1}^{n}(-1)^k\binom{n-1}{k-1} = -(1-1)^{n-1} = 0$$

Thus condition (2) holds. ∎

**Remark 3.3:** For $1 \leq i \leq n$, let $\epsilon_i \in \{1, -1\}$ and $J = \{k: \epsilon_k = 1\}$ then

$$\|\epsilon_1 x_{i_1} + \epsilon_2 x_{i_2} + \cdots + \epsilon_n x_{i_n}\|^2 = \|x_J - x_{I\setminus J}\|^2$$

$$= 2\|x_J\|^2 + 2\|x_{I\setminus J}\|^2 - \|x_J + x_{I\setminus J}\|^2$$

$$= 2\|x_J\|^2 + 2\|x_{I\setminus J}\|^2 - \|x_I\|^2.$$

This gives us a test for identities of the form

$$\sum_{I=\{i_1,i_2,\cdots,i_k\}\subseteq\overline{n}} a_I \sum_{\epsilon_1,\cdots,\epsilon_k \in\{1,-1\}} \|\epsilon_1 x_{i_1} + \epsilon_2 x_{i_2} + \cdots + \epsilon_k x_{i_k}\|^2 = 0$$

(where the inner sum is over all possible choices of signs $\epsilon_1,\cdots,\epsilon_k$) Indeed, the above identity can be transformed to the form

$$0 = \sum_{I\subseteq\overline{n}} a_I \sum_{J\subseteq I} \|x_J - x_{I\setminus J}\|^2 = \sum_{I\subseteq\overline{n}} a_I \sum_{J\subseteq I}\left(2\|x_J\|^2 + 2\|x_{I\setminus J}\|^2 - \|x_I\|^2\right)$$

which has the form that can be verified using Theorem 2.3.

The test for verifying such identities is given in Theorem 3.4. Corollary 3.5 uses this test to prove the parallelopiped law. Namely

$$\sum_{k=1}^{n} \sum_{1 \leq i_1 < i_2 < \cdots < i_k \leq n} (-1)^k 2^{n-k} \sum_{\epsilon_1, \cdots, \epsilon_k \in \{1,-1\}} \|\epsilon_1 x_{i_1} + \epsilon_2 x_{i_2} + \cdots + \epsilon_k x_{i_k}\|^2 = 0 \quad (14)$$

**Theorem 3.4:** Let $\mathcal{H}$ be an inner-product space of dimension at least 2. For each $n \geq 2$, $a_1, \ldots, a_n$ fixed real numbers, and $\mathcal{J} \subseteq \mathcal{P}(\overline{n})$. If $a_I = \prod_{i \in I} a_i$ for $I \subseteq \overline{n}$ then

$$0 = \sum_{I \in \mathcal{J}} a_I \sum_{J \subseteq I} \|x_J - x_{I \setminus J}\|^2 = 0 \quad (15)$$

holds for all $x_1, \ldots, x_n$ in $\mathcal{H}$ if and only if for each $i \in \overline{n}$,

$$0 = \sum \{a_I : I \in \mathcal{J} \text{ and } i \in I\}$$

**Proof:**

As in Remark 3.3, we have

$$0 = \sum_{I \in \mathcal{J}} a_I \sum_{J \subseteq I} \|x_J - x_{I \setminus J}\|^2 = \sum_{I \in \mathcal{J}} a_I \sum_{J \subseteq I} \left(2\|x_J\|^2 + 2\|x_{I \setminus J}\|^2 - \|x_I\|^2\right) \quad (16)$$

For $I$ fixed and $J_1 \subseteq I$ the square norm $\|x_{J_1}\|^2$ occurs (with a factor of 2) twice in the inner sum on the RHS of (16) (once when $J = J_1$ and the other when $J = I \setminus J_1$) while the $-\|x_I\|^2$ occurs $2^{|I|}$ times in the inner sum (once for each $J \subseteq I$) Thus, the RHS of (16) is

$$\sum_{J} 4\|x_J\|^2 \sum_{\{I \in \mathcal{J} : I \supseteq J\}} a_I - \sum_{I \in \mathcal{J}} a_I 2^{|I|} \|x_I\|^2$$

So, the equality in (16) becomes

$$0 = \sum_{J \subseteq \overline{n}} 4 \left(\sum_{\{I \in \mathcal{J} : I \supseteq J\}} a_I\right) \|x_J\|^2 - \sum_{I \subseteq \overline{n}} a_I 2^{|I|} \|x_I\|^2 \quad (17)$$

Using Theorem 2.1 to substitute for $\|x_J\|^2$ and $\|x_I\|^2$ in (17), we get that for each $B$ a subset of $\overline{n}$ of cardinality 2,

$$a_B = \sum_{B \subseteq J \subseteq \overline{n}} 4 \sum_{\{I \in \mathcal{J} : I \supseteq J\}} a_I - \sum_{\{I \in \mathcal{J} : I \supseteq B\}} a_I 2^{|I|}.$$

The first sum is

$$\sum_{B \subseteq J \subseteq I \in \mathcal{J}} 4 a_I = \sum_{B \subseteq I \in \mathcal{J}} 4 a_I \sum_{B \subseteq J \subseteq I} 1.$$

As the sets $J$ satisfying $B \subseteq J \subseteq I$ are in bijective correspondence with the subsets of $I \setminus B$ the number of the former is $2^{|I \setminus B|} = 2^{|I|-2}$, thus

$$a_B = \sum_{B \subseteq I \in \mathcal{J}} a_I 4 \cdot 2^{|I|-2} - \sum_{B \subseteq I \in \mathcal{J}} a_I 2^{|I|} = 0$$

So, condition (1) of Theorem 2.3 is satisfied. Therefore (16) holds if and only if condition (1) of Theorem 2.3 holds. Hence, it suffices to show that in this case condition (2) of Theorem 2.3 is equivalent to the condition in Theorem 3.4. Let $u$ be a unit vector. For $1 \leq i \leq n$ let $x_i = u$ and $x_j = 0$ for $j \neq i$. With this choice of values for the $x_k$,s

For $J \subseteq I$ that $\|x_J - x_{I \setminus J}\|^2 = 0$ if $i \notin I$, and if $i \in I$, we have

$$\|x_J - x_{I\setminus J}\|^2 = \begin{cases} \|x_i\|^2 & \text{if } i \in J \\ \|-x_i\|^2 & \text{if } i \notin J \end{cases}$$

in either case $\|x_J - x_{I\setminus J}\|^2 = \|x_i\|^2 = 1$. Thus

$$\sum_{I \in \mathcal{J}} a_I \sum_{J \subseteq I} \|x_J - x_{I\setminus J}\|^2 = \sum_{\{:I\in\mathcal{J}:\ i\in I\}} \sum_{J \subseteq I} 1 = \sum_{\{:I\in\mathcal{J}:\ i\in I\}} 2^{|I|} a_I$$

proving the desired equivalence. ∎

**Corollary 3.5:** Let $\mathcal{H}$ be an inner-product space of dimension at least 2. For each $n \geq 2$, $x_1, \ldots, x_n$ in $\mathcal{H}$ and real numbers $a_1, \ldots, a_n$, if there are $i \neq j$ such that $a_i = a_j = -1/2$ then

$$0 = \sum_{k=1}^n \sum_{1 \leq i_1 < i_2 < \cdots < i_k \leq n} a_{i_1} a_{i_2} \cdots a_{i_k} \sum_{\epsilon_1, \ldots, \epsilon_k \in \{1, -1\}} \|\epsilon_1 x_{i_1} + \epsilon_2 x_{i_2} + \cdots + \epsilon_k x_{i_k}\|^2 \quad (18)$$

**Proof:**

As in Remark 3.3, the equation in (18) can be rewritten as

$$0 = \sum_{I \in \mathcal{P}(\bar{n})} a_I \sum_{J \subseteq I} \|x_J - x_{I\setminus J}\|^2$$

where, as in Theorem 3.4, $a_I = \prod_{i \in I} a_i$. This is just (15) with $\mathcal{J} = \mathcal{P}(\bar{n})$. Thus, by Theorem 3.4, the identity holds if and only if for each $i \in \bar{n}$, $\sum_{i \in I \subseteq \bar{n}} 2^{|I|} a_I = 0$. For each $k \in \bar{n}$, let $c_k = 2a_k$, then the validity test for our identity can be rewritten as follows:

for each $1 \leq i \leq n$,

$$\sum_{i \in I \subseteq \bar{n}} c_I = c_i \sum_{I \subseteq \bar{n}\setminus\{i\}} c_I = 0.$$

For a nonempty index set $J$, by an easy induction on the cardinality of $J$, we get that

$$\sum_{I \subseteq J} c_I = \prod_{l \in J}(1 + c_l).$$

Thus, the condition of Theorem 3.4 becomes, for each $i \in \bar{n}$

$$2a_i \prod_{i \neq j \in J}(1 + 2a_j) = 0$$

Which always hold if there are two distinct indices $i \neq j$ so that $a_i = a_j = -\frac{1}{2}$ ∎

Note that factoring out $2^n$ from (14), it becomes a special case of Corollary 3.4 where $a_1 = a_2 = \cdots = a_n = -1/2$.

## 4. Sufficient conditions for an inner product

In the previous sections we examined identities that follow from the norm being derived from an inner product. In this section, we show that any such identity (any identity of the form $\sum_{A \subseteq \overline{n}} c_A \|x_A\|^2 = 0$ ) implies that the norm is derived from an inner product. The proof is divided into two lemmas.

**Lemma 4.1:** Let $\mathcal{H}$ be a normed spaced.

1) If some identity of the form $\sum_{B \subseteq \overline{n}} c_B \|x_B\|^2 = 0$ holds in $\mathcal{H}$ with the $c_B \neq 0$ for some $B \neq \emptyset$ then an identity of the form $\sum_{B \subseteq A} a_B \|x_B\|^2 = 0$ with $a_A \neq 0$ and $A \neq \emptyset$ holds in $\mathcal{H}$.
2) If an identity of the form $\sum_{B \subseteq A} a_B \|x_B\|^2 = 0$ with $a_A \neq 0$ and $A \neq \emptyset$ holds in $\mathcal{H}$. Then an identity of the form $\sum_{B \subsetneq A} (-1)^{|A|-|B|} \|x_B\|^2 = 0$ also holds in $\mathcal{H}$ for some $A \neq \emptyset$.

**Proof:**

1) Pick $A$ a maximal element in the collection $\{B \subseteq \overline{n} : c_B \neq 0 \text{ and } B \neq \emptyset\}$ with respect to inclusion. Let $x_i = 0$ for $i \notin A$ and let $x_i$ be arbitrary for $i \in A$. Using this choice of values our identity becomes

$$0 = \sum_{B \subseteq \overline{n}} c_B \|x_{A \cap B}\|^2 = \sum_{B \subseteq A} \|x_B\|^2 \left( \sum_{D \cap A = B} c_D \right)$$

Since $A$ is maximal $D \cap A = A$ ,i.e. $D \supseteq A$, and $c_D \neq 0$ implies that $D = A$ so the coefficient of $\|x_A\|^2$ is $c_A \neq 0$ and $\neq \emptyset$, by our choice of $A$, $A \neq \emptyset$. Thus the identity we obtained above has desired type.
2) Let $N$ be the minimum element of the nonempty set
$\{|A|: A \neq \emptyset \land \text{ there is an identity } \sum_{B \subseteq A} a_B \|x_B\|^2 = 0 \text{ with } a_A \neq 0 \text{ holding in } \mathcal{H}\}$
Pick a set $A$ with cardinality $N$, and an identity $\sum_{B \subseteq A} a_B \|x_B\|^2 = 0$ with $a_A \neq 0$ holding in $\mathcal{H}$.

By our choice of $A$, if $|C| < |A|$ is nonempty and an identity of the form $\sum_{B \subseteq C} a_B \|x_B\|^2 = $ holds in $\mathcal{H}$ then $a_C = 0$ . (*)

We prove that, with the above choice of $A$, $a_B = (-1)^{|A|-|B|} a_A$ for every $\emptyset \neq B \subseteq A$. The proof is by downward induction on $|B|$. For $|B| = |A|$ this is clear since $B = A$. Suppose the claim is true for $k < |B| \leq |A|$ and $k \geq 1$. Let $B_0$ be a subset of $A$ of cardinality $k$. Choose $x_i = 0$ for $i \notin B_0$ and let $x_i$ be arbitrary for $i \in B_0$. Substituting in our identity, we obtain,

$$0 = \sum_{B \subseteq A} a_B \|x_{B \cap B_0}\|^2 = \sum_{C \subseteq B_0} \|x_C\|^2 \left( \sum_{B \cap B_0 = C, B \subseteq A} a_B \right)$$

Since $1 \leq |B_0| < |A|$ we have by (*) that the coefficient of $\|x_{B_0}\|^2$ in the above expression is zero, so $\sum_{B \cap B_0 = B_0, B \subseteq A} a_B = 0$. Therefore, by the induction hypothesis,

$$a_{B_0} = - \sum_{A \supseteq B \supsetneq B_0} a_B = - \sum_{\supseteq B \supsetneq B_0} (-1)^{|A|-|B|} a_A.$$

For each $m > k$ the number of m-subsets $B \subseteq A$ that contain $B_0$ is $\binom{|A|-k}{m-k}$ so

$$a_{B_0} = -a_A \sum_{m=k+1}^{|A|} \binom{|A|-k}{m-k} (-1)^{|A|-m}$$

$$= -a_A \left( \sum_{m=k}^{|A|} \binom{|A|-k}{m-k}(-1)^{|A|-m} - (-1)^{|A|-k} \right)$$

$$= -a_A \left( (1-1)^{|A|-k} - (-1)^{|A|-k} \right) = (-1)^{|A|-k} a_A$$

Establishing the result for subsets of cardinality $k$. Thus for our choice of $A$ the above identity is $a_A \sum_{\emptyset \neq B \subsetneq A}(-1)^{|A|-|B|}\|x_B\|^2 = \sum_{B \subsetneq A}(-1)^{|A|-|B|}\|x_B\|^2 = 0$, and we obtain the desired conclusion by dividing by $a_A \neq 0$. ∎

**Lemma 4.2:** Let $\mathcal{H}$ be a normed space. If for some $A \neq \emptyset$ the identity

$$\sum_{B \subsetneq A}(-1)^{|A|-|B|}\|x_B\|^2 = 0$$

holds in $\mathcal{H}$, then $\mathcal{H}$ is an inner product space.

**Proof:**

WLOG assume that $A = \bar{n}$ for some $n \in \mathbb{N}$. If $n = 3$ our hypothesis (after renaming variables) coincides with the identity (1) from which we obtain identity (2) by replacing $z$ by $-y$. The existence of inner product then follows from the Jordan and von Neumann result.

If $n > 3$, Let $x, y \in \mathcal{H}$ and $t \in \mathbb{R}$ be arbitrary. The substitution

$$x_1 \to x \text{ and } x_k \to ty \text{ for } 2 \leq k \leq n$$

in the identity in the statement of the lemma, yields

$$\sum_{j=0}^{n-1} \binom{n-1}{j}(-1)^{n-j}\|x + jty\|^2 = 0. \tag{18}$$

Indeed, for $1 \in B$ and $|B| = j+1$, the above substitution transforms $(-1)^{|A|-|B|}\|x_B\|^2$ into $(-1)^{n-j-1}\|x + jty\|^2$.

Thus,

$$\sum_{1 \in B \subseteq A}(-1)^{|A|-|B|}\|x_B\|^2 = \sum_{j=0}^{n-1}(-1)^{n-j-1}\binom{n-1}{j}\|x + jty\|^2 \tag{19}$$

For $1 \notin B$, we have $x_B = |B|ty$. Thus, since there are $\binom{n-1}{j}$ j-subsets of $A$ that don't contain 1

$$\sum_{1 \notin B \subsetneq A}(-1)^{|A|-|B|}\|x_B\|^2 = t^2\|y\|^2 \sum_{j=0}^{n-1}(-1)^{n-j}j\binom{n-1}{j}$$

Using the identity $j\binom{n-1}{j} = (n-1)\binom{n-2}{j-1}$, we get that the last sum is

$$(n-1)t^2\|y\|^2 \sum_{j=1}^{n-1}(-1)^{n-j}\binom{n-2}{j-1} = -(n-1)t^2\|y\|^2(1-1)^{n-2} = 0$$

Thus,

$$0 = -\sum_{B \subseteq A}(-1)^{|A|-|B|}\|x_B\|^2 = \sum_{j=0}^{n-1}(-1)^{n-j}\binom{n-1}{j}\|x + jty\|^2$$

Fix $x$ and $y$. The function $g(t) = \|x + ty\|^2$ is continuous since if $K$ is a number such that $|s|, |t| \leq K$, then $|g(t) - g(s)|$ is bounded above by

$$|\|x + t y\| - \|x + s y\|| 2(\|x\| + K \|y\|) \leq 2|s - t| \|y\|(\|x\| + K \|y\|)$$

The sum

$$\sum_{j=0}^{n-1} (-1)^{n-j} \binom{n-1}{j} g(r + js)$$

is just

$$\sum_{j=0}^{n-1} (-1)^{n-j} \binom{n-1}{j} \|x + (r + js)y\|^2 = \sum_{j=0}^{n-1} (-1)^{n-j} \binom{n-1}{j} \|x_1 + j s y\|^2$$

where $x_1 = x + ry$ which is zero by (18) so, by Lemma 1.1, $g$ is a polynomial. Since for $k > 2$,

$$\lim_{t \to \infty} \frac{g(t)}{t^k} \leq \lim_{t \to \infty} \frac{1}{t^k}(\|x\| + t\|y\|)^2 = \lim_{t \to \infty} \frac{1}{t^{k-2}} \left(\frac{1}{t^2}\|x\| + \|y\|\right)^2 = 0$$

the degree of $g$ is at most 2. Thus,

$$g(t) = A + Bt + Ct^2.$$

we have, $A = g(0) = \|x\|^2$ and $B = \lim_{t \to \infty} \frac{g(t)}{t^2} = \lim_{t \to \infty} \|\frac{x}{t} + y\|^2 = \|y\|^2$.

Thus,

$$\|x + t y\|^2 = g(t) = \|x\|^2 + B t + \|y\|^2$$

and

$$\|x + y\|^2 + \|x - y\|^2 = g(1) + g(-1) = 2\|x\|^2 + 2\|y\|^2.$$

Since $x, y \in \mathcal{H}$ were arbitrary, by Jordan and von Neumann result, the norm is induced by an inner product. ∎

The following result follows directly from Lemmas 4.1, and Lemma 4.2.

**Theorem 4.3:** Let $\mathcal{H}$ be a normed space. If an identity of the form $\sum_{B \subseteq \bar{n}} c_B \|x_B\|^2 = 0$ holds in $\mathcal{H}$ with the $c_B \neq 0$ for some $B \neq \emptyset$, then the norm is given by an inner product.

**Corollary 4.5:** Let $\mathcal{H}$ be a normed space. The norm on $\mathcal{H}$ is given by an inner product if and only if any of any of the identities in Lemmas 3.1 or 3.2 or Corollary 3.5 hold.

**Proof:**

If $\mathcal{H}$ is an inner-product space, then the identities hold by the Lemmas in which they occur. If, on the other hand, any of these identities hold then by Theorem 4.3, $\mathcal{H}$ is an inner-product space. ∎

**Acknowledgment:**

The author would like to thank the Deanship of scientific research at King Faisal university for their continuous support during this research work.

**References**

[1] Mirosslaw Adamek, Characterization of Inner Product Spaces by Strongly Schur-Convex Functions, Results in Mathematics 75, 72 (2020) https://doi.org/10.1007/s00025-020-01197-1

[2] Dan Amir, Characterizations of inner product spaces, Birkhäuser Verlag 1986.

[3] G. Chelidze, Inner product spaces and minimal values of functionals, Journal of Mathematical Analysis and Applications, Volume 298, Issue 1, 1 October 2004, Pages 106-113

[4] F. Dadipour and M. S. Moslehian, A characterization of inner product spaces related to the p-angular distance, Journal of Mathematical Analysis and Applications, Volume 371, Issue 2, 15 November 2010, Pages 677-681.

[5] M. Fréchet, Sur la definition axiomatique d'une classe d'espaces vectorials distanciés applicables vectoriellement sur l'espaces de Hilbert, Annals of Mathematics Second Series, vol.36, no.3, (1935) pp. 70-78.

[6] Istrâţescu, V.I Inner Product Structures: Theory and Applications, D. Reidel Publishing Company, 1987.

[7] P. Jordan and J. von Neumann, On inner products in linear, metric spaces, Annals of Mathematics Second Series, vol.36 no.3, (1935) pp. 719-730.

[8] José Mendoza and Tijani Pakhrou, On some characterizations of inner product spaces, Journal of Mathematical Analysis and Applications Volume 282, Issue 1, 1 June 2003, Pages 369-382.

[9] B. Reznick, Banach spaces which satisfy linear identities. Pacific Jour. Math. vol 74, (1978) pp. 221-233.